\documentclass[a4paper, 11pt] {article} 
\usepackage{amssymb}
\usepackage{url}
\usepackage{graphicx}

\begin{document}

\title{Effect of population size in a Prey-Predator model}
\author{Fabien Campillo\thanks{EPI Modemic INRA/INRIA, SupAgro, 2 place Viala, 34060 Montpellier Cedex 2, France} \and Claude Lobry\footnotemark[1]
}

\maketitle

\begin{abstract}
We consider a stochastic version of the basic predator-prey differential equation model. The model, which contains a parameter $\omega$ which represents the number of individuals for one unit of prey - If $x$ denotes the quantity of prey in the differential equation model $x=1$ means that there are $\omega$ individuals in the discontinuous one - is derived from the classical birth and death process. It is shown by the mean of simulations and explained by a mathematical analysis based on results in {\em singular perturbation theory} (the so called  theory of {\em Canards}) that qualitative properties of the model like persistence or extinction are dramatically sensitive to $\omega$. For instance, in our example, if $\omega = 10^7$ we have extinction and if $\omega = 10^8$ we have persistence. This means that we must be very cautious when we use continuous variables in place of jump processes in dynamic population modeling even when we use stochastic differential equations in place of deterministic ones.

\end{abstract}

\section{Introduction}
\noindent Consider the standard prey-predator model :

\begin{equation}\label{syst1}
\Sigma\quad\quad \left  \{
 \begin{array}{lcl}
 \displaystyle \frac{dx}{dt}&=&f(x) -\mu(x)y  \\[8pt]
\displaystyle \frac{dy}{dt}&=& (c\,\mu(x) - \delta)y
 
 \end{array} 
\right .
\end{equation}
where $x$ stands for the concentration of preys and $y$ for the concentration of predators. It is well known that this kind of modeling with differential equations is valid only if {\em one unity} of $x$ (or $y$) represents a {\em large number } of prey (or predator) {\em individuals}. 
On the other hand, when the number of individual is too small, everybody agree that one must switch to some kind of individually based modeling of stochastic nature.\\\\
What means {\em large} is generally not specified but it is widely admitted that around $10^3$ the law of large numbers begins to do its job and that figures like $10^6$ are completely safe if one wants to use continuous variables and differential equations.\\\\
The objective of this paper is to show that the threshold of $10^3$ is not always acceptable and that, in some circumstances, even $10^6$ cannot be considered as secure when we deduce biological consequences, like persistence, from the behavior of a model with continuous variables. For that purpose we propose
a stochastic model, where the dynamic of the prey is governed by a birth and death process while, for mathematical simplicity, we keep the predator variable as a continuous one. The development will make clear that this simplification does not affect the conclusions of the paper. The proposed model is such that the dynamic of the process is locally approximated (when the number of preys is large) by a differential system which is precisely of predator-prey type like (\ref{syst1}). We agree that, in many respects, our model is biologically questionable but our objective is not to contribute to biological understanding of prey-predator relationship but just to point out some mathematical phenomenon which is likely to be present in many models and which might be responsible for erroneous interpretations.\\\\
The first section is devoted to the presentation of the stochastic model, the second to the presentation of some surprising simulations, the third to the analysis of the differential system that governs the dynamics of the mean of the stochastic process and the forth to the explanations of the surprising aspects of the simulations. The last two sections are devoted to methodological  and bibliographical  comments.\\\\
From the mathematical point of view the material and results presented here are classical. The paper is intended principally for non mathematically oriented readers who are not necessarily aware of these questions. We tried to avoid all mathematical technicalities and for this purpose we made an important use of results from computer simulations. All the references to existing literature related to these questions are rejected to the last two sections.

 
\section{The model.}
\noindent The variable $\omega\,x(t)$ is an integer which is the number of preys at time $t$. This variable performs the following birth and death (actually here ``death'' means ``capture'' by a predator) process.
\begin{itemize}
\item At any time, the epoch $\tau$ of the next event (birth or death) is a random variable $Z$ which follows an exponential law of parameter :

\begin{equation}\label{lambda}
\lambda = \frac{\omega}{\varepsilon}(f(x)+\mu(x)y)
\end{equation}
\item At the epoch $\tau$ we have one birth with probability $\frac{f(x)}{f(x)+\mu(x)y}$ or one death with the complementary probability :
\begin{equation} \label{proba}
 \begin{array}{lcr}
\displaystyle P(\omega \,x(\tau^+) = \omega \,x(\tau^-) + 1) &= &\displaystyle\frac{f(x(\tau^-))}{f(x(\tau^-))+\mu(x(\tau^-))y(\tau^-)} \\[12pt]
\displaystyle P(\omega \,x(\tau^+) = \omega \,x(\tau^-) - 1) &=&\displaystyle \frac{\mu(x(\tau^-))y(\tau^-)}{f(x(\tau^-))+\mu(x(\tau^-))y(\tau^-)}
 \end{array} 
 \end{equation}
\end{itemize}
The variable $y$ is a continuous variable which evolves according to :
\begin{equation}\label{vary}
y(t+dt) = y(t) - dt\,m\, y(t) + \varepsilon \{ \mathrm{number \;of\; captures\; during} \; [t, t+dt] \}
\end{equation}
Thus the predator dynamics is an exponential decay associated to a growth proportional to the number of prey disappearing during the elapsed time. The parameter $\varepsilon$ accounts for different time scale for the prey and the predator dynamics.\\\\
Assume that $dt = 10^{-4}$, $\omega = 10^9$, $\varepsilon = 10^{-1}$ and $f(x)+\mu(x)y$ is of the order of unity. Then, during elapsed 
time $dt$ the number of events (death or birth) is of the order of $\lambda \; dt = \frac{\omega}{\varepsilon}(f(x)+\mu(x)y)dt \approx \frac{10^5}{\varepsilon}\approx 10^6 $. This is a bit lengthy to simulate (at least with a desk computer) but, due to that great number of events, the process defined by (\ref{lambda}), (\ref{proba}), (\ref{vary}) is accurately approximated on the interval $[t, t+dt]$ by the {diffusion process} (see appendix A for a derivation) :

\begin{equation}\label{systd}
 \left  \{
 \begin{array}{lcl}
 \displaystyle x(t+dt) &=&x(t) +\frac{dt}{\varepsilon}[f(x(t)) -\mu(x(t))y(t) ]-\sigma_x W_t  \\[8pt]
 \displaystyle y(t+dt)&=& y(t) + dt [(\mu(x(t))-m)y(t)]+ \sigma_y W_t
 
 \end{array} 
\right .
\end{equation}
where $ W_{1dt},W_{2dt},W_{3dt},......$ is a sequence of independent Gaussian variables with mean $0$ and standard deviation $1$ with :

$$\sigma_x = \sqrt{ \frac{4 dt}{\omega \varepsilon}}\;\; \sqrt{\frac{f(x(t)\mu(x(t))y(t)}{f(x(t))+\mu(x(t))y(t)}}$$

$$\sigma_y = \sqrt{\frac{ dt \varepsilon}{\omega}}\;\; \sqrt{\frac{f(x(t)\mu(x(t))y(t)}{f(x(t))+\mu(x(t))y(t)}}$$
This diffusion process is not a good approximation of the jump process when $x$ is small. For an accurate description one must switch  to the jump process for small values of $x$ but, since this is not our point here, we restrict us to the consideration of the
stochastic (with continuous variables) diffusion like process :
\begin{equation}\label{syst2}
 \left  \{
 \begin{array}{lcl}
\;\; \mathrm{if  }\;\; x(t)& \leq \frac{1}{\omega}& \mathrm{then }\;\;\; x(t+dt) = 0 \;\;\;\mathrm{ else} \\
 \displaystyle x(t+dt) &=&x(t) +\frac{dt}{\varepsilon}[f(x(t)) -\mu(x(t))y(t) ]+\sigma_x W_t  \\[8pt]
 \displaystyle y(t+dt)&=& y(t) + dt [(\mu(x(t))-m)y(t)]+ \sigma_y W_t
 \end{array} 
\right .
\end{equation}
The first line in (\ref{syst2}) states that when the number of prey is smaller than $1$ it has to be $0$ . It is necessary to specify this because now the variable $x$ is continuous in the diffusion process but we want to keep the meaning of $x$ as a {\em number of individuals} . Thus, $\frac{1}{\omega}$ must be an absorbing barrier for (\ref{syst2}). For $x \geq \frac{1}{\omega}$ one sees that the recurrence equation for the mean of $x(t)$ and $y(t)$ is approximated by :
\begin{equation}\label{syst3}
 \left  \{
 \begin{array}{lcl}
E[x(t+dt)] &= &E[x(t)] + \frac{dt}{\varepsilon} [f(x(t))-\mu(x(t))y(t)]\\
E[y(t+dt)] &= &E[y(t)] + dt[(\mu(x(t)-m)y(t)]
 \end{array} 
\right .
\end{equation}
which is the Euler scheme for the differential system :
\begin{equation}\label{syst4}
\left  \{
 \begin{array}{lcl}
 \displaystyle \frac{dx}{dt}&=&\frac{1}{\varepsilon}[f(x) -\mu(x)y ] \\[8pt]
 \displaystyle \frac{dy}{dt}&=& (\mu(x) - m)y
 \end{array} 
\right .
\end{equation}
Thus, to conclude this paragraph, we have constructed a diffusion-like model defined by equations (\ref{syst2}). This model depends on a parameter $\omega$. This model has the following properties :
\begin{itemize}
\item Since the model is derived from the birth and death process $\omega\,x$ is interpreted as the number of individuals for $x$ units of preys.
\item The standard deviation is proportional to $\sqrt{\frac{1}{\omega}}$ : the biggest is $\omega$ the more ``deterministic'' is the process.
\item The diffusion process is degenerate (i.e. the dimension of the random noise is not $2$ but $1$). This is due to the fact that only $x$ is considered as a discrete variable, not $y$.
\item When $\omega\,x$ is large (greater than $10^3$) the dynamic of the model is accurately approximated (at least for small durations) by the classical deterministic differential prey-predator model (\ref{syst4})  which is the same as system  (\ref{syst1}) with $ c = \varepsilon, \; \varepsilon \, m = \delta $ after a change of time units. 
\end{itemize}
We shall first simulate this system and then explain the observed simulations.
\section{Simulations}
\begin{figure}[h,t,b]   \centering
   \includegraphics[width=12cm]{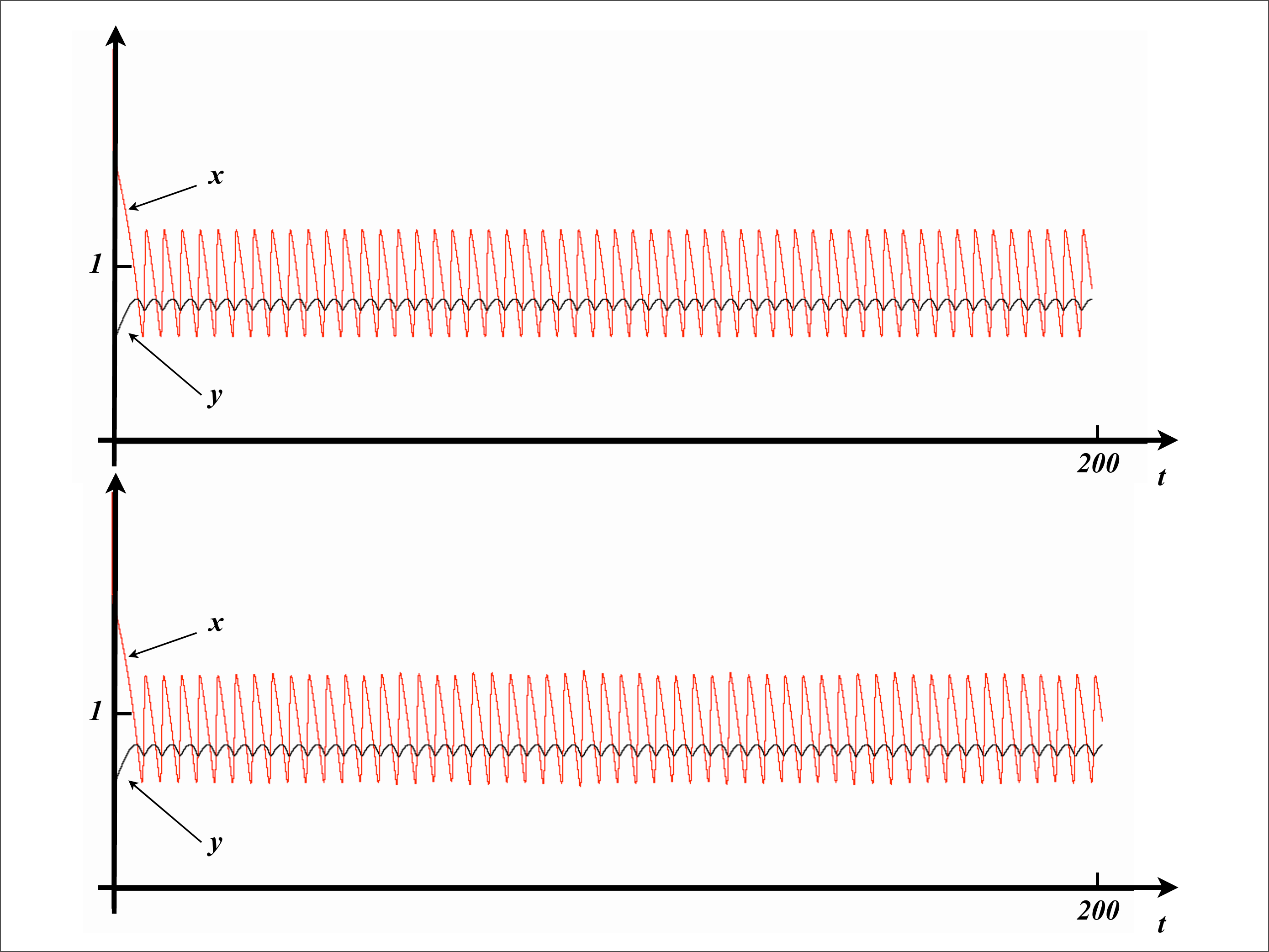} 
   \caption{Above : $\omega = 10^{12}$ ; below $\omega = 10^{10}$}
   \label{fig-1}
\end{figure}

\begin{figure}[h,t,b]   \centering
   \includegraphics[width=12cm]{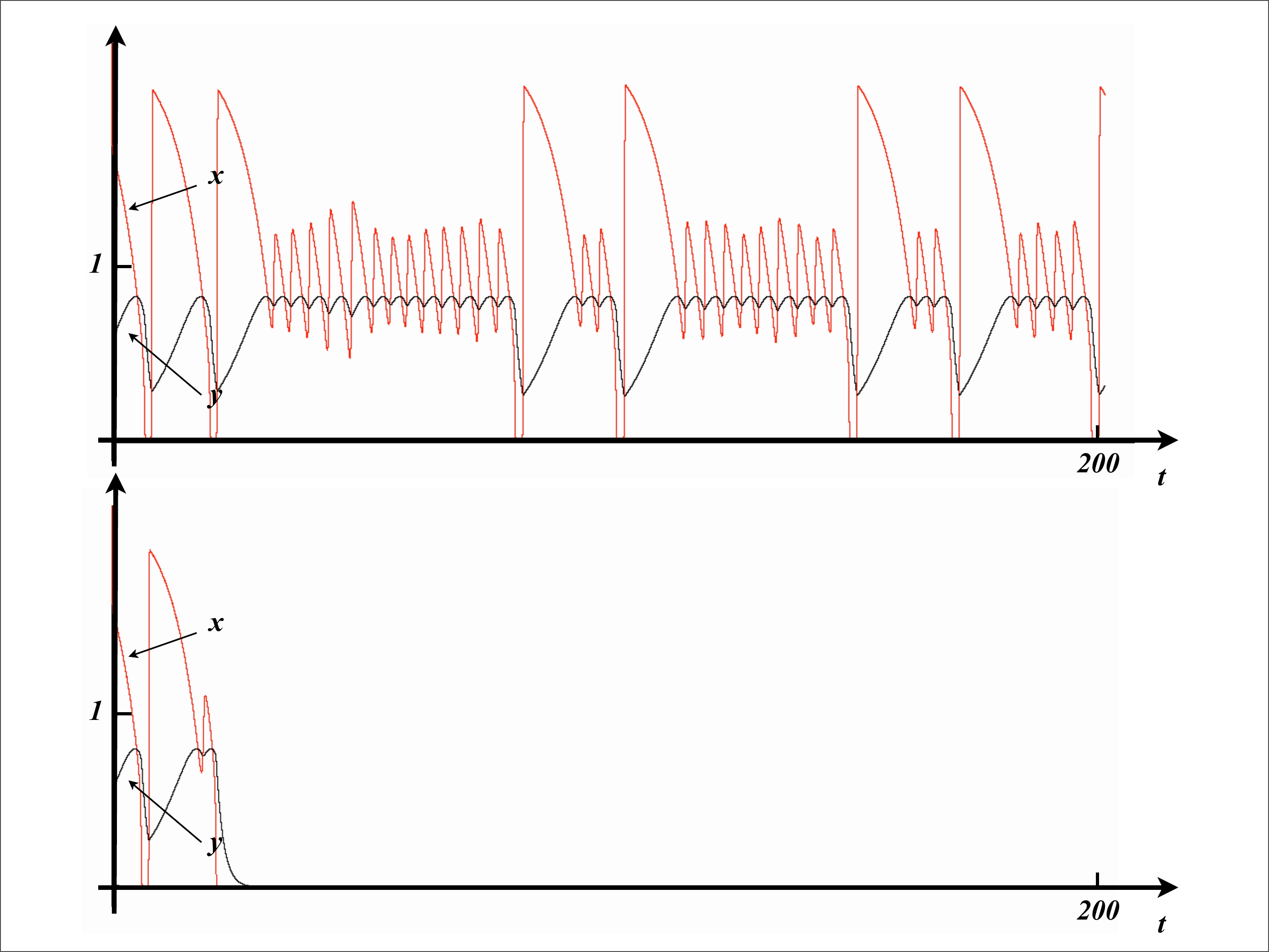} 
   \caption{Above : $\omega = 10^8$ ; below $\omega = 10^6$}
   \label{fig-2}
\end{figure}
\begin{figure}[h,t,b]   \centering
   \includegraphics[width=12cm]{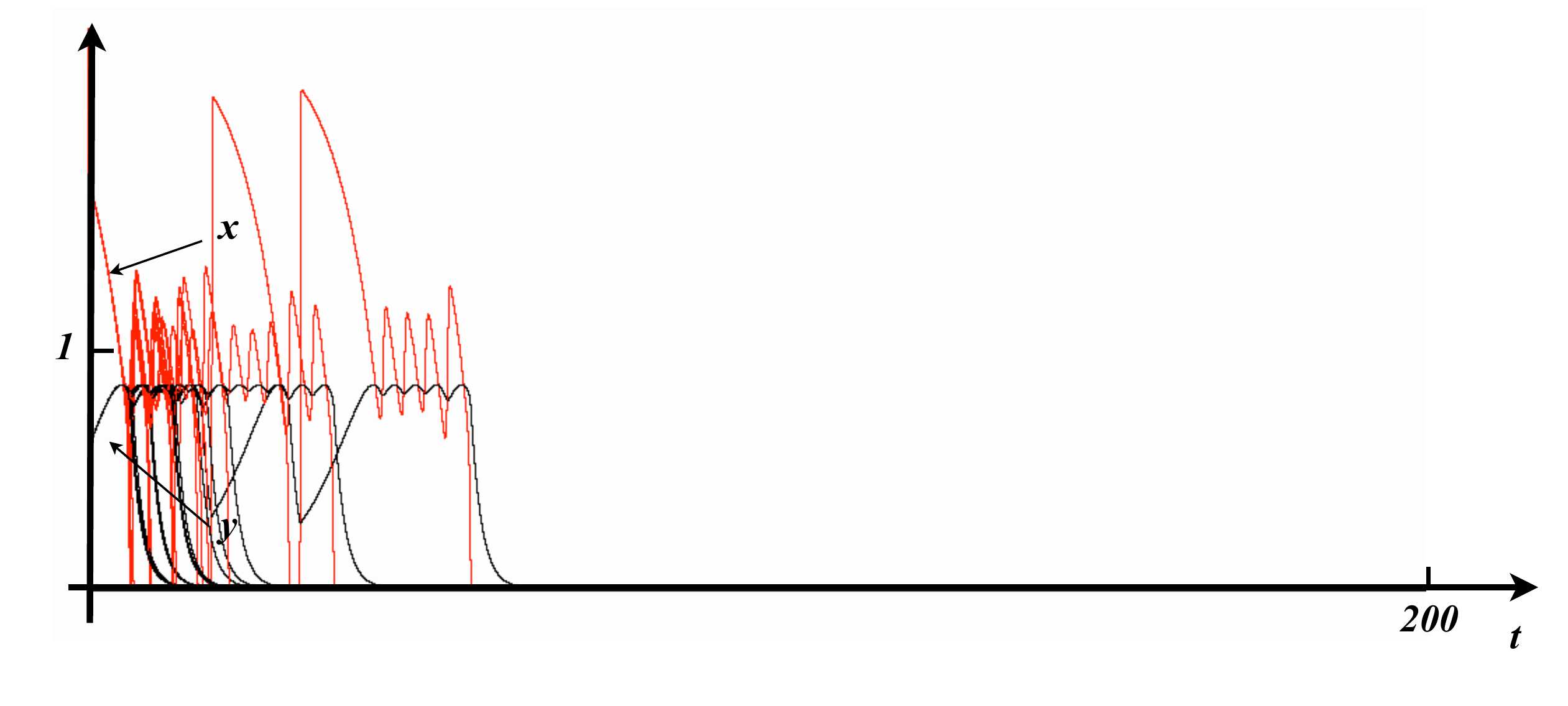} 
   \caption{Twenty runs with  $\omega = 10^6$ }
   \label{fig-3}
\end{figure}
\noindent In this section we fix $f$, $\mu$, $\varepsilon$ and $m$ as : 
\begin{itemize}
\item $ f(x) = \displaystyle \frac{1}{2} x\;(2-x)$
\item $\mu (x) = \displaystyle \frac{x}{0.4+x}$
\item $\varepsilon = 0.02$
\item $ m = 0.6645$
\end{itemize} 
On Fig.\ref{fig-1} one sees one run of the process (\ref{syst2}). The duration is $200$ and $\omega$ is fixed at $10^{12}$ (above) and $10^{10}$ (below). One sees regular oscillations for the population of prey (in red) and predator (in black). We do not see any difference between the two records. These regular oscillations are those predicted by the deterministic prey-predator model. Since the value of $x$ during the oscillations is around $1$ which corresponds to such a great number of individuals we are definitely not surprised that the continuous deterministic system is a good approximation.\\\\
But, on Fig.\ref{fig-2} we observe a dramatic change with $\omega = 10^8$ {\em which is still a big figure}. We observe a {\em mixed mode oscillation} with a random successions of large and small oscillations which could not be produced by a deterministic  two dimensional system. With $\omega = 10^6$ we observe an extinction of the two populations which is confirmed on Fig.\ref{fig-3} where we observe that none of $20$ runs for  $\omega = 10^6$ is persistent at time $T = 200$.\\\\
Let us denote by $T$ the time of extinction for the predator (defined as the time where $y(t)$ reaches the value $\frac{1}{\omega}$).
Let us say that there is extinction when the time of extinction is smaller than $1000$. On Tab. \ref{table-1} we have the empirical probabilities of extinction with respect to $\omega$ (computed on 1000 runs) and mean and standard deviation of $T$ computed on trajectories ending with extinction for $t < 1000$ . 
\begin{table}[h,t,b]   
 \begin{center}
      \begin{tabular}{|r|r|r|r|}
   \hline 
    $\omega$&$E[T]$&$\sigma(T)$&$P(T\leq 1000$) \\
    \hline
    $10^5$&30.46&6.75&1\\
    \hline
     $10^6$&39.02&11.30&1\\
      \hline
      $2.0\;10^6$&47.74&19.62&1\\
    \hline
    $4.0\;10^6$&79.05&51.79&1\\
    \hline
     $6.0\;10^6$&143.54&121.42&0.999\\
    \hline
     $8.0\;10^6$&259.76&222.64&0.983\\
    \hline
    $9.0\;10^6$&311.70&247.58&0.964\\
    \hline
 $10^7$&554.21&319.94&0.867\\
    \hline
    $1.1\; 10^7$&555.17&351.75&0.741\\
    \hline
    $1.2\; 10^7$&681.31&324.12&0.649\\
    \hline
   $1.3\; 10^7$&745.83&321.95&0.481\\
    \hline
        $1.4\; 10^7$&815.46&296.26&0.384\\
    \hline
$1.5\; 10^7$&867.54&273.60&0.255\\
    \hline
    $1.6\; 10^7$&906.10&238.55&0.182\\
\hline
$1.7\; 10^7$&928.68&221.50&0.120\\
\hline
$1.8\; 10^7$&964.05&143.82&0.072\\
\hline
$1.9\; 10^7$&975.48&110.09&0.059\\
\hline
$2.0\;10^7$& $> 1000$ &&0\\
    \hline
   \end{tabular}
   \end{center}
   \caption{Empirical probabilities of extinction according to $\omega$}  \label{table-1}
\end{table}
We can see that the transition is very sharp from extinction with probability one ($\omega = 4.0\,10^6$) to non extinction with probability one ($\omega = 2.0\,^7$). It seems surprising
that with about $\omega = 2.0\,10^7$ prey-individuals the system is definitely (say up to 1000 units of time) safe and definitely unsafe for $4.0\,10^6$ which is still a big figure. This is a problem since in most case, in population dynamics models, we have poor information on the actual size of a population. We come back later on this issue. Notice also that the standard deviation of $T$ is very large for small values of $\omega$ which makes predictions very imprecise.

\section{The dynamics of the continuous deterministic model}
\begin{figure}[h,t,b]   \centering
   \includegraphics[width=8cm]{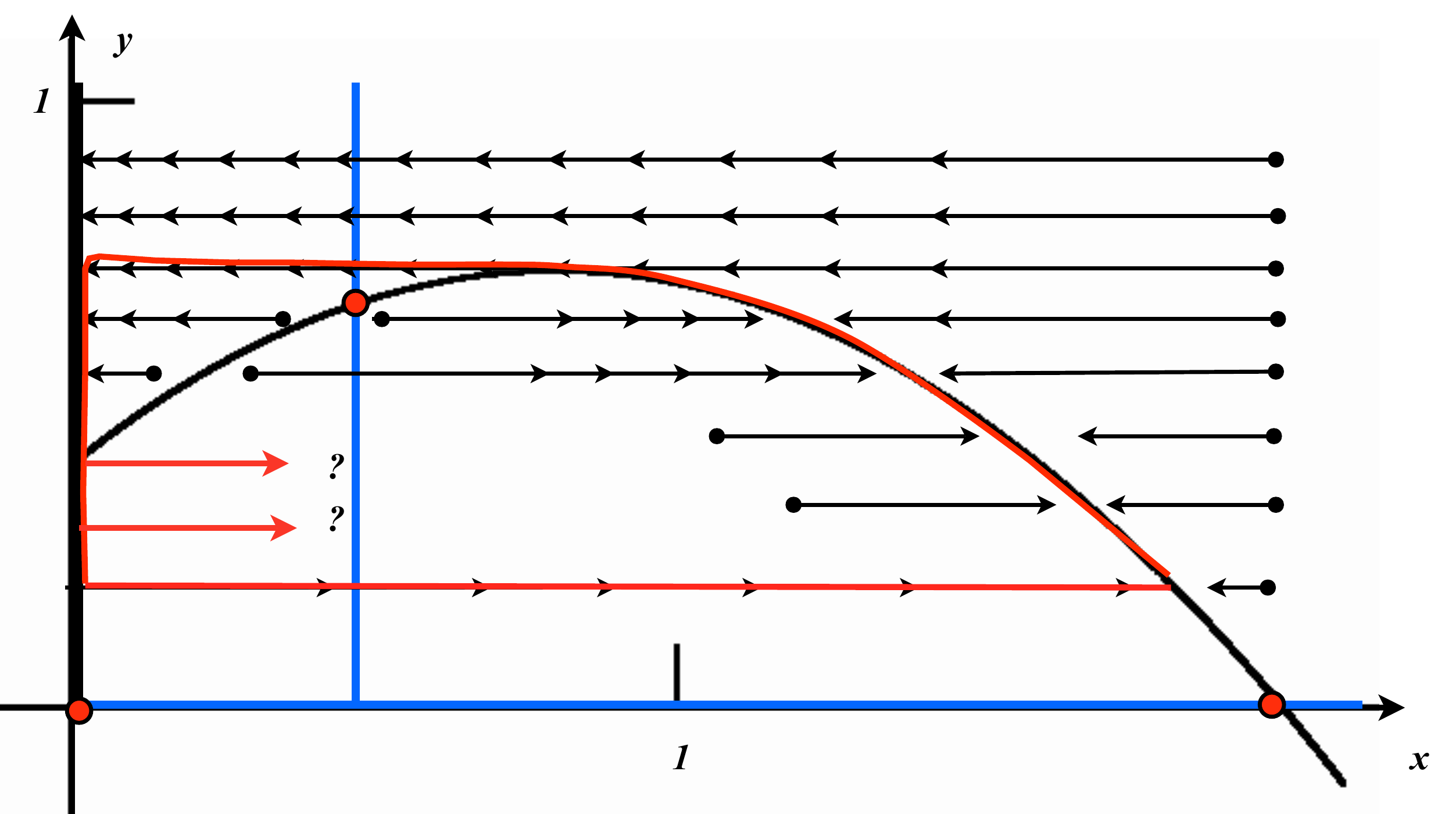} 
   \caption{}
   \label{schem1}
\end{figure}

\begin{figure}[h,t,b]   \centering
   \includegraphics[width=8cm]{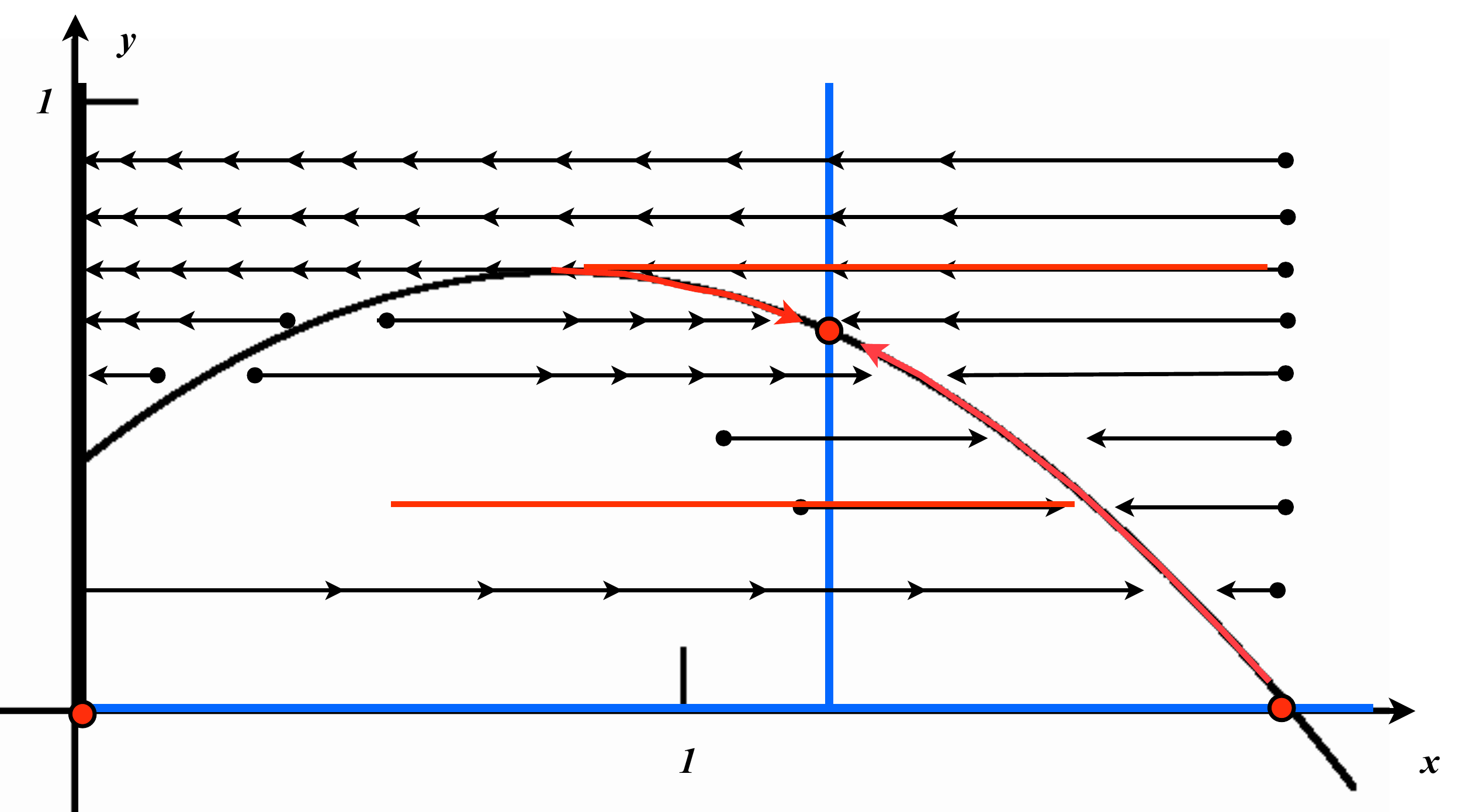} 
   \caption{Schematic representation of solutions of (\ref{syst4}) }
   \label{schem2}
\end{figure}
\noindent In this section we describe the dynamics of the deterministic model (\ref{syst4}) which approximate the evolution of the mean of the diffusion model (\ref{syst3}). All the material in this section is classical and known as the theory of ``canards" (see the section ``literature comments" for more details).\\\\
The first step in the understanding of a planar system like our is to draw the two nullclines (sometimes called ``zero growth isoclines''),  that is the sets defined by  :
\begin{itemize}
\item The nullcline of the prey : $\{ (x,y) : \frac{1}{\varepsilon}\;[f(x) -\mu(x)y ]= 0\}$
\item The nullcline of the predator :$\{ (x,y) : (\mu(x)-m)y = 0\}$
\end{itemize}
In our simulations the parameter $\varepsilon$ is small ($0.02$) and, by the way, except when the quantity 
$$[f(x) -\mu(x)y ] $$
is small, of the order of $\varepsilon$, the right member in the first equation in (\ref{syst4}) is large compared to the second one. This means that the vector velocity of (\ref{syst4}) is {\em almost horizontal}. From this it follows that, a first approximation the solutions of our system is shown by the hand-drawn schemes on Fig.\ref{schem1} and Fig.\ref{schem2} : Outside of the parabola and the $y$ axe which is the nullcline of the prey 
the trajectories are taken as horizontal.
\begin{itemize}
\item On Fig.\ref{schem1} one sees that the nullcline of the predator (in blue) is on the left of the maximum of the  nullcline of the prey (the black parabola curve). Along the nullcline of the prey the motion is down-up on the right of the vertical blue nullcline and up-down on the left.
From this we see that there is a tendency for the trajectories to join the $y$-axe on its attractive part (above the black curve), to follow it in the up-down direction and there is some indeterminacy to where it will leave it after having crossed the the parabola. From this scheme we suspect the existence of a periodic  limit-cycle cycle which, actually, can be proven to be the case.
\item On Fig.\ref{schem2} the situation is somewhat easier to understand. The blue vertical null-cline of the predator being on the right the motion along the parabola converges to a limit point which apparently is a stable attracting point for all initial conditions.
\item Notice that an attracting equilibrium and an attracting limit-cycle are qualitatively different picture and that the transition between the two cases occurs when $m$ crosses the value $0.666....$ (when the blue line crosses the parabola at its maximum).
\end{itemize}
\begin{figure}[h,t,b]   \centering
   \includegraphics[width=12cm]{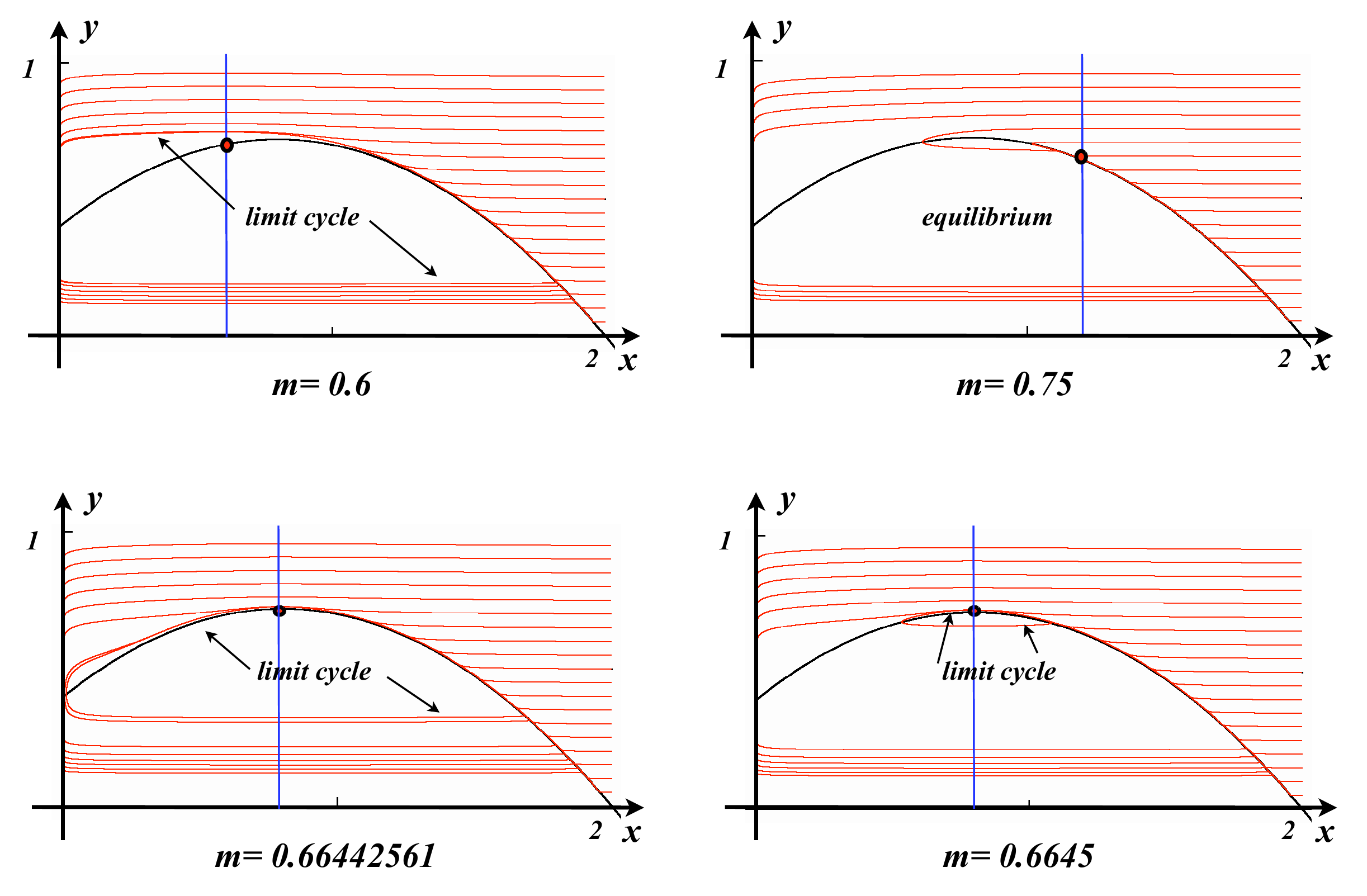} 
   \caption{Phase portrait of system (\ref{syst4}) for different values of $m$}
   \label{Bif1}
\end{figure}
Let us now comment on Fig.\ref{Bif1}. The pictures are not hand-drawn schemes but actual simulations with $\varepsilon = 0.02$ ; we observe the great similarity with the schemes.
\begin{itemize}
\item $m = 0.6$ : One large limit cycle (the direction of the motion is counter clockwise). Trajectories above the limit cycle are of two kinds : some hit the limit cycle and then follow it and the others hit the vertical axe, then they follow it up-down and reappear below, run  left to right hit the parabola and then join the limit cycle. Actually ``true'' trajectories never meet but, due to the limit of our drawing, they seem to meet. All trajectories follow for a while the $y$ axe and then $x(t)$ is potentially small.
\item $m = 0.75$ : We have an attracting equilibrium. Some trajectories go directly to the equilibrium, some other follow the $y$ axe.
\item $m = 0.6645$ : We have a small periodic limit cycle circling around the unstable equilibrium which is very close to the periodic orbit. Along the periodic orbit we stay very far from the $y$ axe (and, by the way, $x(t)$ is never small) but one sees that near the
unstable equilibrium, a very small perturbation leads to a trajectory which follows the $y$ axe and $x(t)$ can become small.
\item $m = 0.66442561$ : In this case we have a limit cycle which is of intermediate size between ``large'' (follows the $y$ axe for a while and small (remains far fro the $y$ axe) ; it just hits the $y$-axe. The point is that it needs very sharp values for $m$ (8 digits in our case) to obtain this intermediate cycle called a ``canard cycle''. See in the section ``comments'' some informations about the mathematical theory of ``canards''.
\end{itemize}
All along this description we said that $x(t)$ is potentially small when the trajectory follows the $y$ axe. {\em But how small ?} A simple
way to enlarge what is going on along this axe is to plot, not the point $(x(t),y(t))$, but $ (\xi (t), y(t))$ with :
$$ \xi (t) = \varepsilon \mathrm{ln}(x(t))$$
\begin{figure}[h,t,b]   \centering
   \includegraphics[width=10cm]{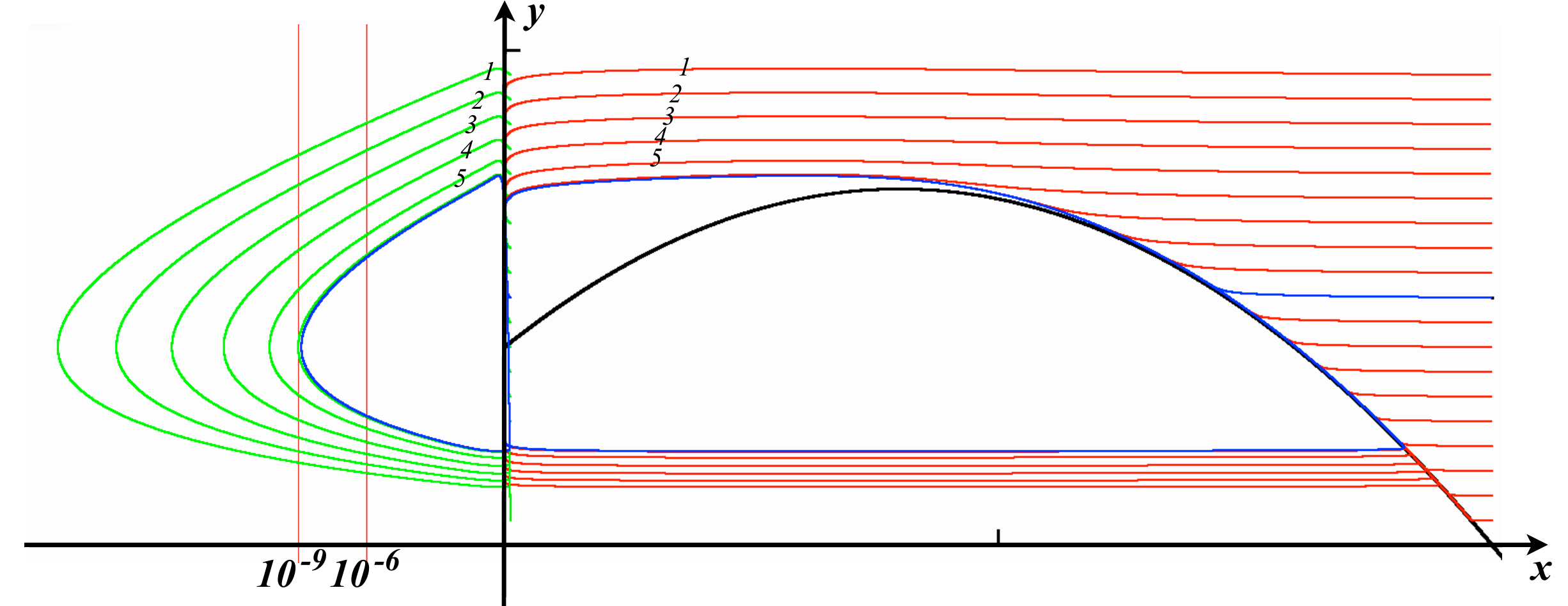} 
   \caption{$(x,y)$ and $(\xi,y)$ variables on the same axes : $m=0.6$}
   \label{log1}
\end{figure}
\begin{figure}[h,t,b]   \centering
   \includegraphics[width=10cm]{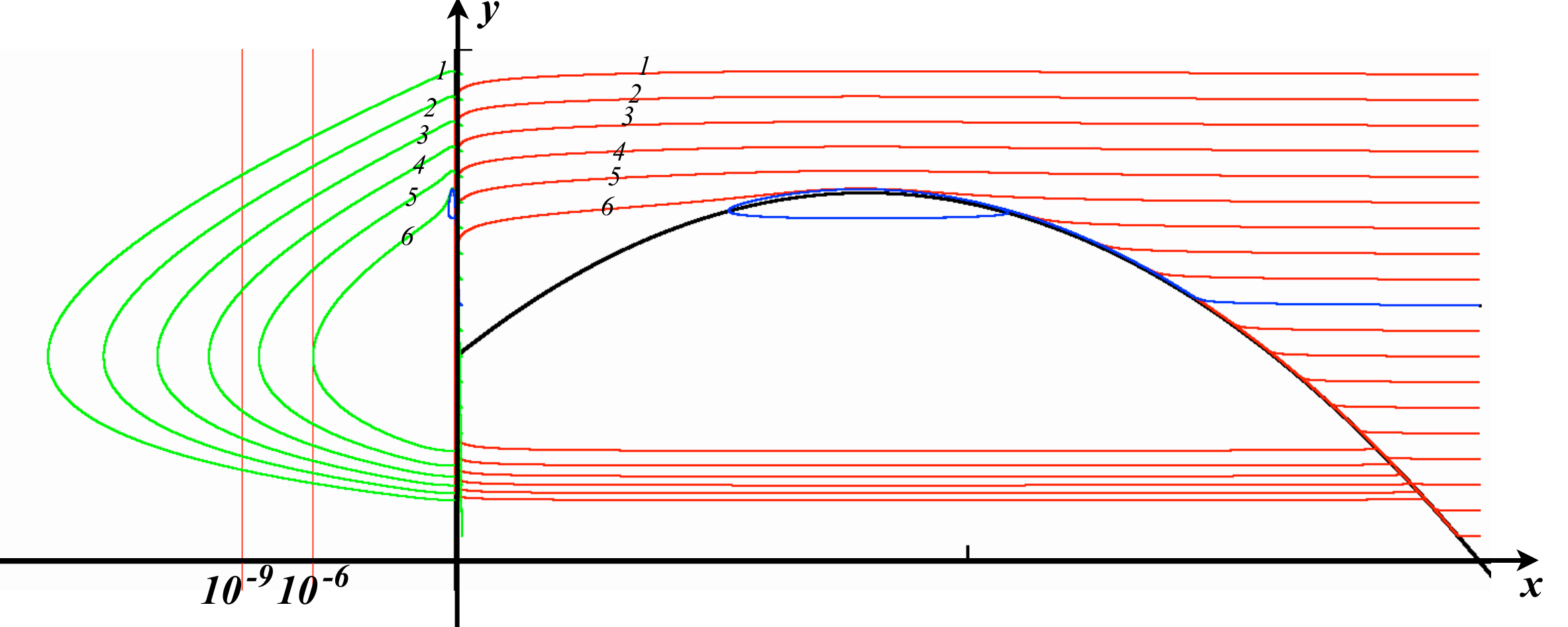} 
   \caption{$(x,y)$ and $(\xi,y)$ variables on the same axes : $m=0.6445$}
   \label{log2}
\end{figure}This is done on Fig.\ref{log1} and Fig.\ref{log2}. We represent the $(x,y)$ and and the $(\xi , y)$ trajectories in the same system of axes ;  $(x,y)$ trajectories are in red, $(\xi , y)$ are in green and both limit cycles in the two systems of representation in blue. The two vertical red lines correspond to $x=10^{-9}$ and $10^{-6}$. There are $19$ trajectories starting from $(2, 0.5\pm k\, 0.05)\;k = 0, 1, ...,9$.\\\\
Let us compare the two simulations.
\begin{itemize}
\item \textbf{Fig.\ref{log1}}. We look at the ``large'' limit cycle in the $(\xi,y)$ variables and we see that the minimum of $\xi$ corresponds to $x = 10^{-9}$ ; for the trajectory labeled $1$ the minimum is about $10^{-17}$.  These incredibly small values are easily explained in Appendix B.
\item \textbf{Fig.\ref{log2}}. The ``small'' limit cycle is almost not visible in the $(\xi,y)$ variables. In both the $(\xi,y)$ and the $(x,y)$ variables the trajectories labeled $1$ to $4$ look very similar. But, in Fig.\ref{log2}, trajectory $5$ remains above $10^{-9}$  which is not the case in Fig.\ref{log1} and  trajectory $6$ does not exist in Fig.\ref{log1} . 
\end{itemize}
The main difference between the case $m= 0.6$ and the case $m = 0.6645$ is that, in the first case, every trajectory is such that the minimum of $x$ is smaller than $10^{-9}$ unlike in the second case where there are two set of trajectories : Those that start above trajectory $6$ for which the minimum will be smaller than $10^{-6}$ before reaching the limit cycle and the others for which $x(t)$ remains greater than $10^{-6}$. Notice that this trajectory $6$ is in some places very close to the limit cycle.\\\\
The observed differences between $m = 0.6$ and $ m = 0.6645$ are not specific of these values. In particular the same behavior with two type of trajectories separated by a sharp transition is true for all values of $m$ between $m=0.66442561$ and $m=0.6666....$. This behavior is summarized by the description of the ``safety funnel''  shown in Fig.\ref{funnel} by the green arrow and that we explain now. Assume that for some reason we do not accept to pursue a trajectory such that the min of $x(t)$ is smaller than $\alpha = 10^{-k}$ (it may be because we think that the size of the population is to low in order to survive or because we want to switch to a different - stochastic - model). The form $10^{-k}$ is by no mean essential for $\alpha$, it is just to emphasize that $\alpha$ is small. It exists a unique $y_0$ such
that the solution issued from $(\infty, y_0)$ (in practice $2$ is a good infinite), which we call the ``$\alpha$-safety trajectory'', is such that $x(t)$ first decreases and attains a first local minimum equal to $10^{-k}$.  This is the red trajectory on the scheme of  Fig.\ref{funnel}.  This trajectory, when $x(t) \approx \mu^{-1}(m)$  will be very close to the limit cycle (the blue trajectory) . We call $\rho (\varepsilon, k)$ the distance between the two curves ; this can be evaluated from the value of $\varepsilon$ and $k$. The ``safety funnel'' is defined by the parts of red and blue curves on the right of the vertical $x = \mu^{-1}(m)$. If a trajectory which enters the funnel is perturbed, as long as it remains in the funnel, the (future) minimum of x will remain greater the $10^{-k}$. If not, there is a danger to reach values smaller than $10^{-k}$.
\begin{figure}[h,t,b]   \centering
   \includegraphics[width=10cm]{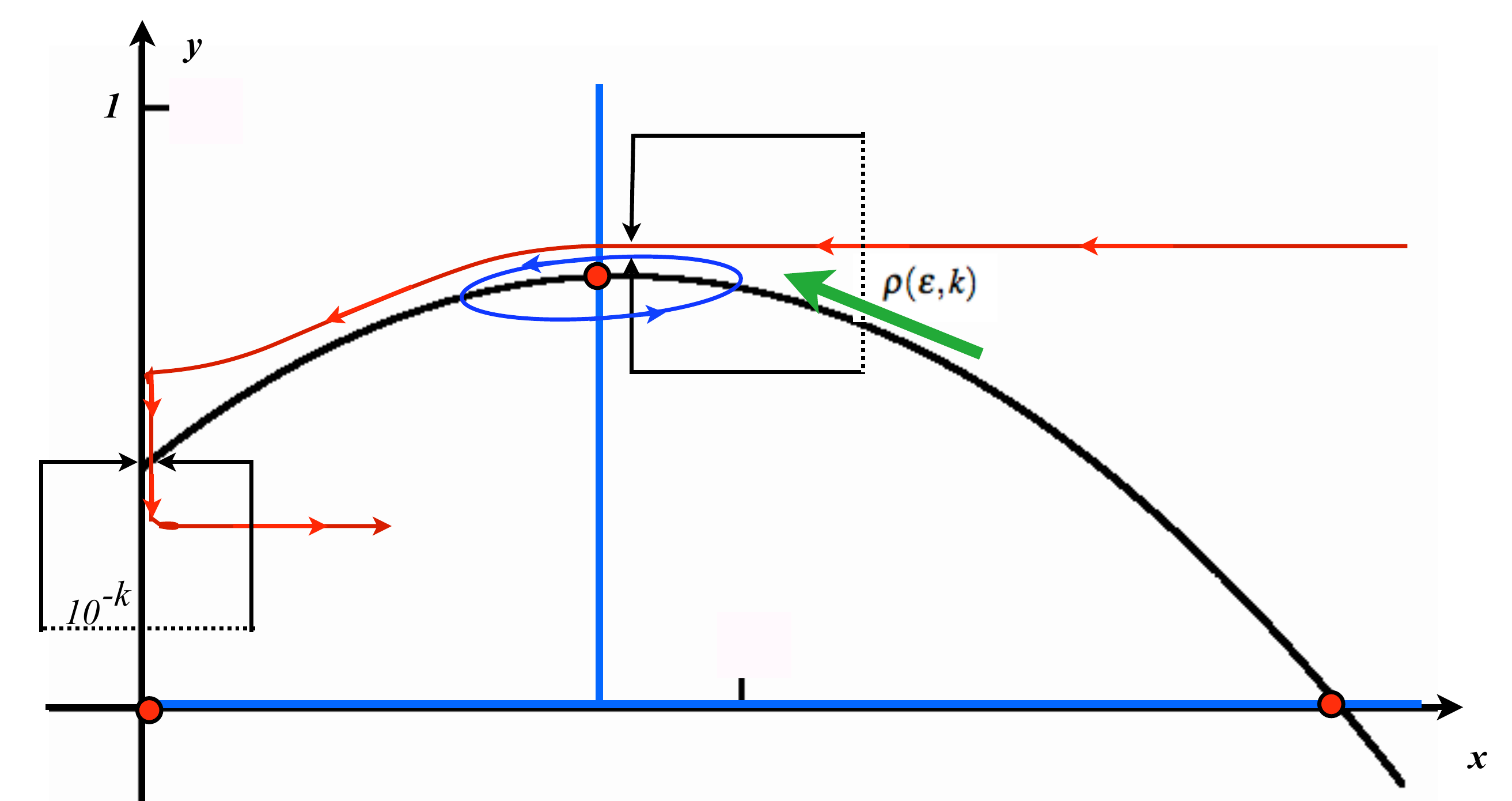} 
   \caption{The ``safety funnel''}
   \label{funnel}
\end{figure}

\section{The diffusion process in the variables $(x,y)$ and $(\xi, y)$}
 \begin{table}[h,t,b]   
\begin{center}
      \begin{tabular}{|r||r|r|r|r|}
   \hline 
    $\omega$ &   $10^9$ &  $10^8$ &  $10^7$ &  $10^6$ \\ 
    \hline
    $\rho(\varepsilon,k)$&   $1.2\;10^{-3}$&  $9.0\;10^{-5}$& $5.5\;10^{-5}$&  $5.3\;10^{-5}$\\
          \hline
          $\sigma_x \sqrt{dt}$&   $4.2\;10^{-5}$ &  $1.4\;10^{-5}$& $4.2\;10^{-5}$&  $1.4\;10^{-4}$\\
          \hline
$\sigma_y \sqrt{dt}$&   $4.9\;10^{-9}$&  $1.4\;10^{-7}$& $4.9\;10^{-7}$&  $1.4\;10^{-6}$\\
          \hline
      \end{tabular}
       \end{center}
      \caption{Width of the funnel and corresponding $\sigma_x$ and $\sigma_y$. }  \label{table-2}
      \end{table}

\begin{figure}[h,t,b]   \centering
   \includegraphics[width=10cm]{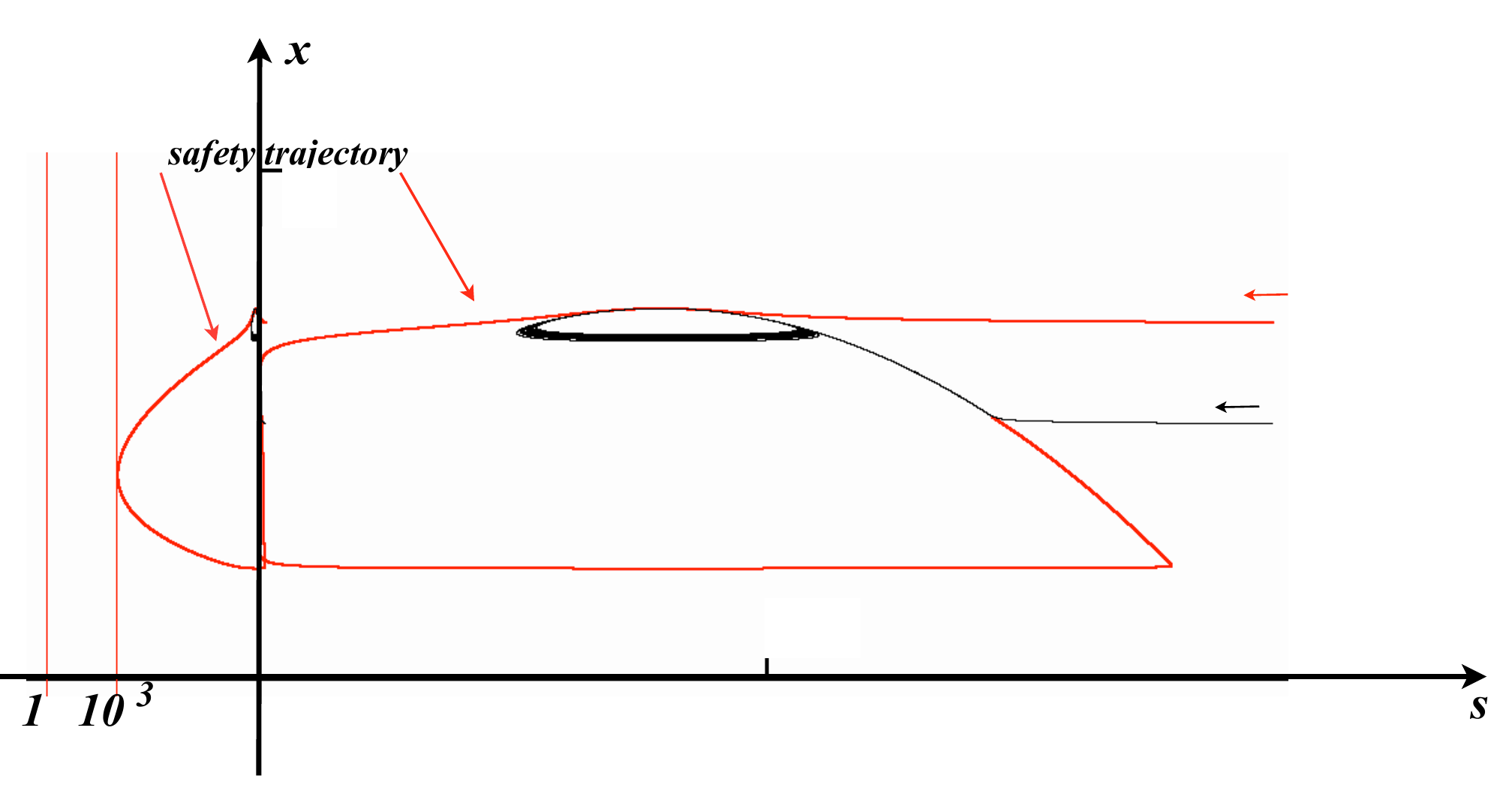} 
   \caption{$\;\;\omega = 10^9$}
   \label{10-9}
\end{figure}
\begin{figure}[h,t,b]   \centering
   \includegraphics[width=10cm]{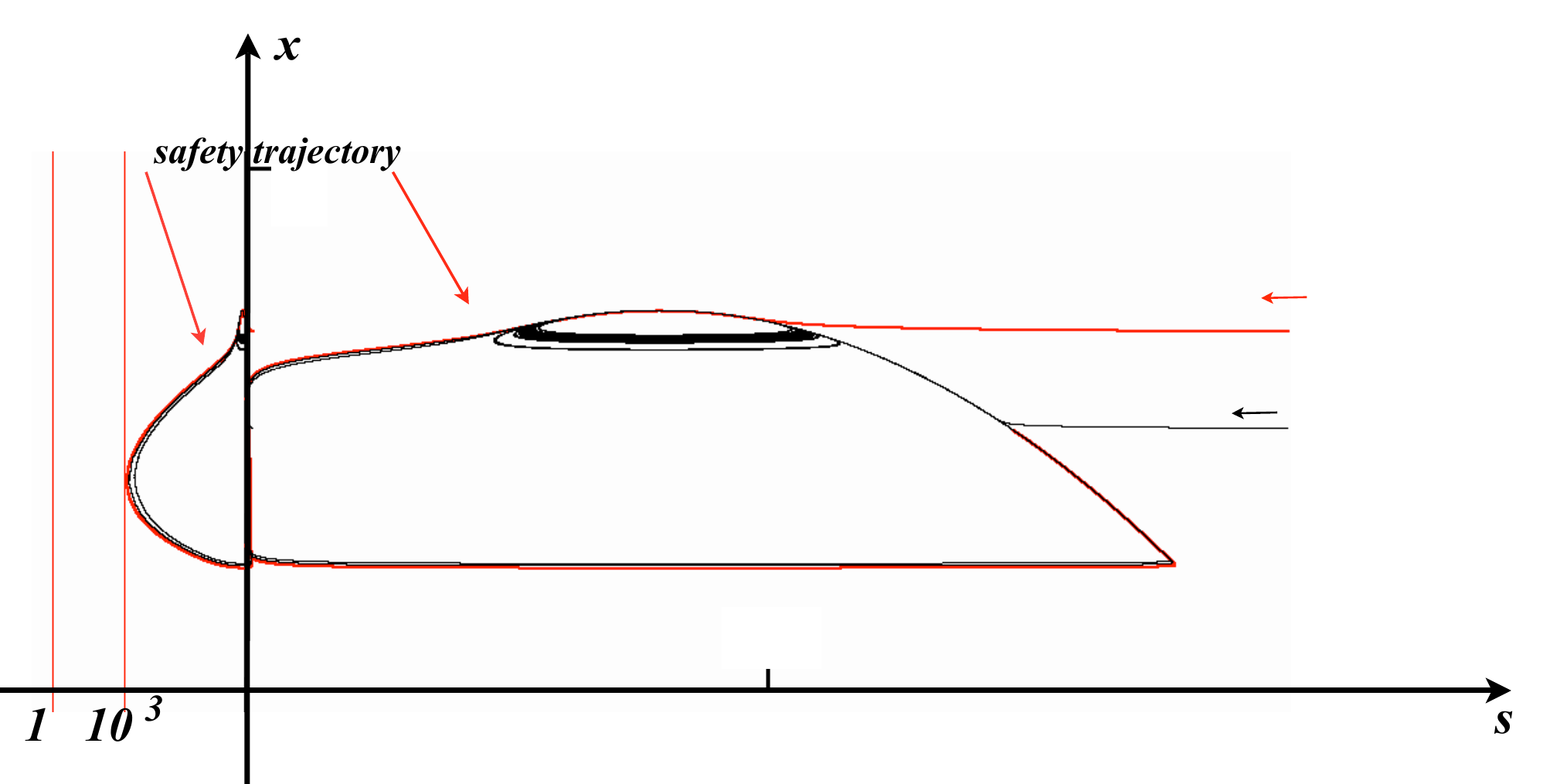} 
   \caption{$\omega = 10^8$}
   \label{10-8}
\end{figure}
\begin{figure}[h,t,b]   \centering
   \includegraphics[width=10cm]{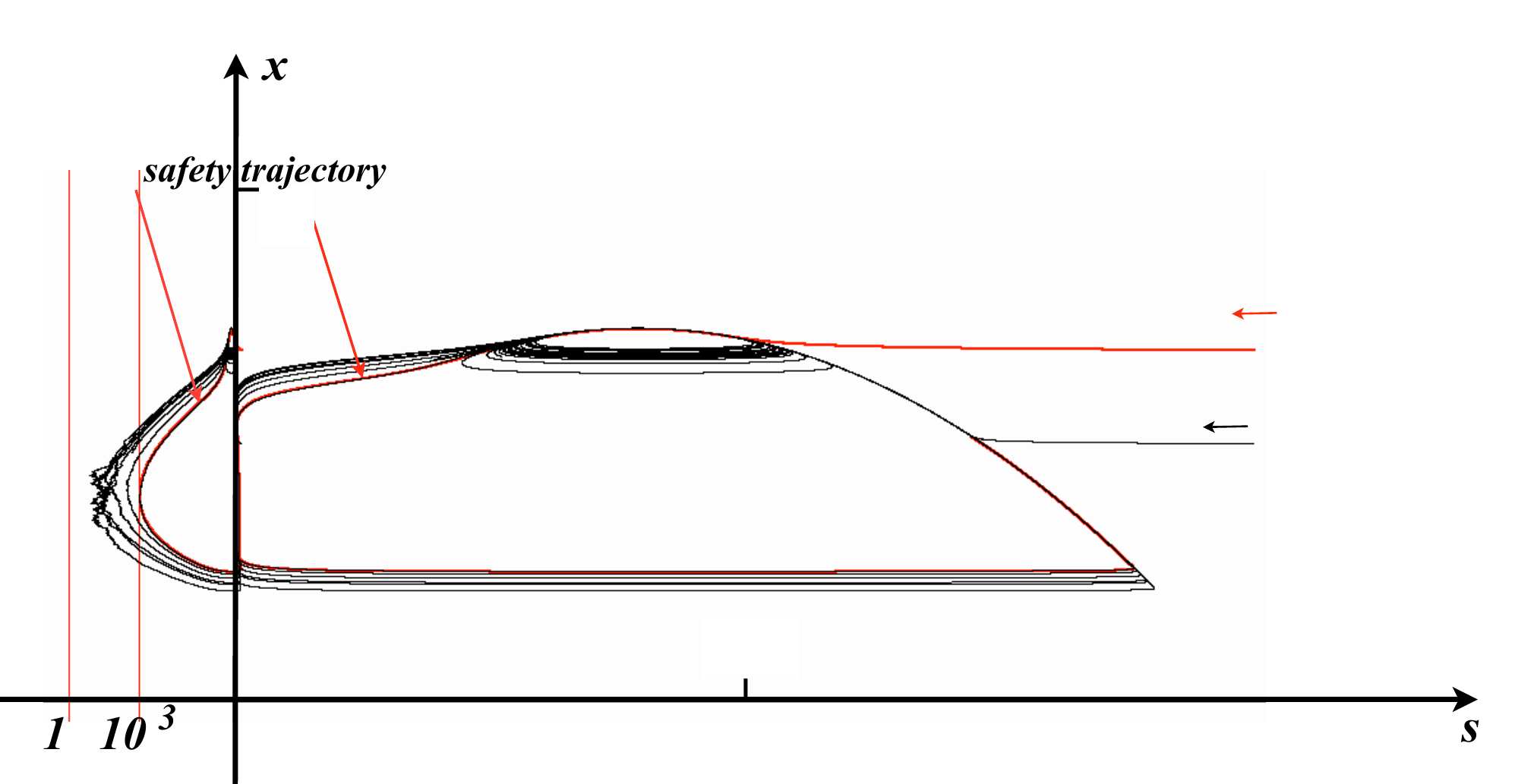} 
   \caption{$\omega = 10^7$}
   \label{10-7}
\end{figure}
\begin{figure}[h,t,b]   \centering
   \includegraphics[width=10cm]{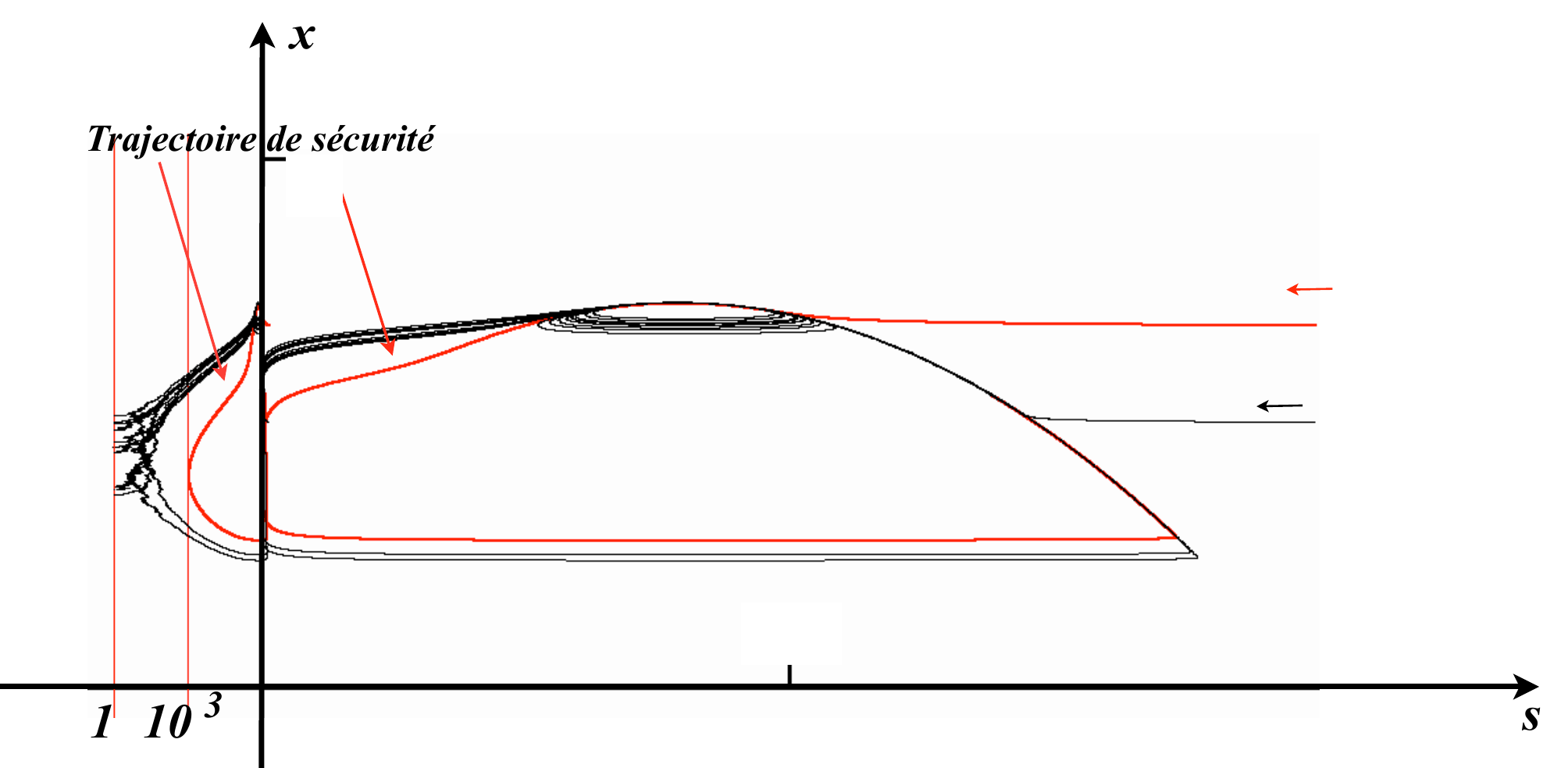} 
   \caption{$\omega = 10^6$}
   \label{10-6}
\end{figure}

\noindent On the four simulations shown on Fig.\ref{10-9} to \ref{10-6} we have performed $10$ runs of  $20$ time units duration of the process (\ref{syst2}) starting from $(-2,0.5)$. The results are presented in both $(x,y)$ and $(\xi,y)$ variables (black trajectories). In the same variables we have simulated from system (\ref{syst4}) the ``safety trajectory'' corresponding to $1000$ individuals, that is to say :
The ``safety trajectory" was obtained by dichotomy and the observed width $\rho(\varepsilon, k)$ of the funnel is given on the table \ref{table-2} with corresponding rough evaluations of $\sigma_x$ and $\sigma_y$ around the funnel. 
\begin{itemize}
\item \textbf{Fig.\ref{10-9}}: All the runs are  widely below the ``safety trajectory''.
\item \textbf{Fig.\ref{10-8}}: All the runs are below the ``safety trajectory'' but we observe that some runs are close to it.
\item \textbf{Fig.\ref{10-7}}: All the runs are above the ``safety trajectory'' and some are closse to the vertical line  corresponding to $10$ individual.
\item \textbf{Fig.\ref{10-7}}: All the runs are widely above the ``safety trajectory'' and  reach ultimately the vertical line corresponding to $1$ individual.
\end{itemize}
We observe that when $\omega$ decreases the strength of the randomness increases and at the same time the width of the funnel decreases. These opposite trends are responsible for the sharp transition from extinction to persistence as $\omega$ grows from $4.0\,10^6$ to $2.0\,10^7$.
\section{Methodological comments}
\subsection{About the question of size of populations.}
\noindent We are used to the fact that continuous differential models works rather well in fluid dynamics and chemical kinetics despite the ultimate discrete nature of fluids. We know that this efficiency is related to the very large number of atoms in the process. Von Foerster, Lotka, Volterra and others popularized the formalism of chemical kinetics in the domain of population dynamics ; they were certainly aware of the limits of such an approach but, in the absence of computers and with a far less developed probability theory, it was a way to progress.\\\\ 
Now, thanks to computers and probability theory, we have good models for small populations. Unfortunately these models are still expensive in term of computer time and deterministic or diffusion models (stochastic differential equations) with continuous variable are still unavoidable. In a diffusion model the size of the population considered is directly related to the strength (standard deviation) of the random term.\\\\
{\em The example of prey-predator interaction presented here shows that the qualitative behavior of such models may depend strongly on the size of the population even when it is very large.}
\subsection{About the generality of the example}
\noindent The deterministic prey-predator model (\ref{syst4}) approximate the dynamics for $E[x]$ and $E[y]$ of the birth and death model defined by (\ref{lambda}), (\ref{proba}), (\ref{vary}) and its diffusion approximation (\ref{syst2}). This model (\ref{syst4}) is the very classical deterministic prey-predator model which is proposed in every text book as a first improvement of the Lotka-Volterra model. The separation of time scales for prey and predator dynamics introduced by the presence of the parameter $\varepsilon$ in the first equation has the following classical explanation. Using a change of time unit (\ref{syst4}) rewrites : 
$$
\left  \{
 \begin{array}{lcl}
 \displaystyle \frac{dx}{d\tau}&=&[f(x) -\mu(x)y ] \\[8pt]
 \displaystyle \frac{dy}{d\tau}&=& \varepsilon (\mu(x) - \delta)y\;\;\;;\;\;\; \varepsilon m = \delta 
 
 \end{array} 
\right .
$$
If we use the same mass unit for $x$ and $y$ then $\varepsilon$ is a yield factor. A yield factor like $0.02$ is acceptable in ecology (one needs $50$ kg. of dry grass to get $1$ kg. of cow). For bigger $\varepsilon$ like $0.1$ the sharp transition that we presented is still present but less spectacular.\\\\
As previously said we admit that our birth and death model is questionable with respect to its biological signification. There are certainly many different models for individual behavior with the same deterministic equation approximating the mean of the process. Since our point relies on the diffusion approximation for such models our conclusions are valid as long as such approximation is correct. In the case of birth and death processes it works provided that the number of individuals is greater than $10^3$-$10^4$ which is our case. For more elaborated models at the individual scale (for instance physiologically structured preys) this point remains to be considered.
\subsection{About the existence of ``canard '' solutions in the model.}
\noindent Let us say two words about  ``canard solutions ''. In a system with two time scales like :
$$
\left  \{
 \begin{array}{lcl}
 \displaystyle \frac{dx}{dt}&=&\frac{1}{\varepsilon}[f(x,y)] \\[8pt]
 \displaystyle \frac{dy}{dt}&=& g(x,y) 
 \end{array} 
\right .
$$
consider the curve $\Gamma$ defined by the equation $f(x,y) = 0$ ; this curve split in two regions :
\begin{itemize}
\item the attracting one made of points such that, in the neighborhood, the vector field converges to $\Gamma$,
\item the repelling one made of points where,  in the neighborhood, the vector field diverges from $\Gamma$,
\end{itemize}
separated by equilibria. A {\em ``canard ''} solution is a solution of the differential system which follows, for some duration, the attracting part of $\Gamma$ at a distance of the order of $\varepsilon$ and, after that, {\em follows also} the repelling part at a distance of order $\varepsilon$. Some `` canards '' are {\em robust } which means that they persist under small changes in the model, others are not.\\\\
The presence of a ``safety funnel'' like the one described in section $3$ is related to the presence of two ``canard solutions'' in (\ref{syst4}).
\begin{itemize}
\item The solution $t\rightarrow (x(t) = 0, y(t) = y(0) e^{-mt})$ which corresponds to the absence of prey,
\item a solution following the cubic from the right to the left, which has no analytic expression but which existence can be proved by continuity arguments.
\end{itemize}
The first ``canard'' is robust but the second is not. This is the reason why, the sharp transition between $4.0\,10^6$ and $2.0\,10^7$ individuals 
occurs for a rather short interval of values of the parameter $m$. As a consequence, to some extend, our example is exceptional, not ``generic''. This will be the case in most two dimensional systems, but this do not invalid our point since robust ``canard'' (different from
trivial ``canard'' corresponding to the absence of some population) are generically present for dimension $3$ and more. 
\begin{figure}[h,t,b]   \centering
   \includegraphics[width=10cm]{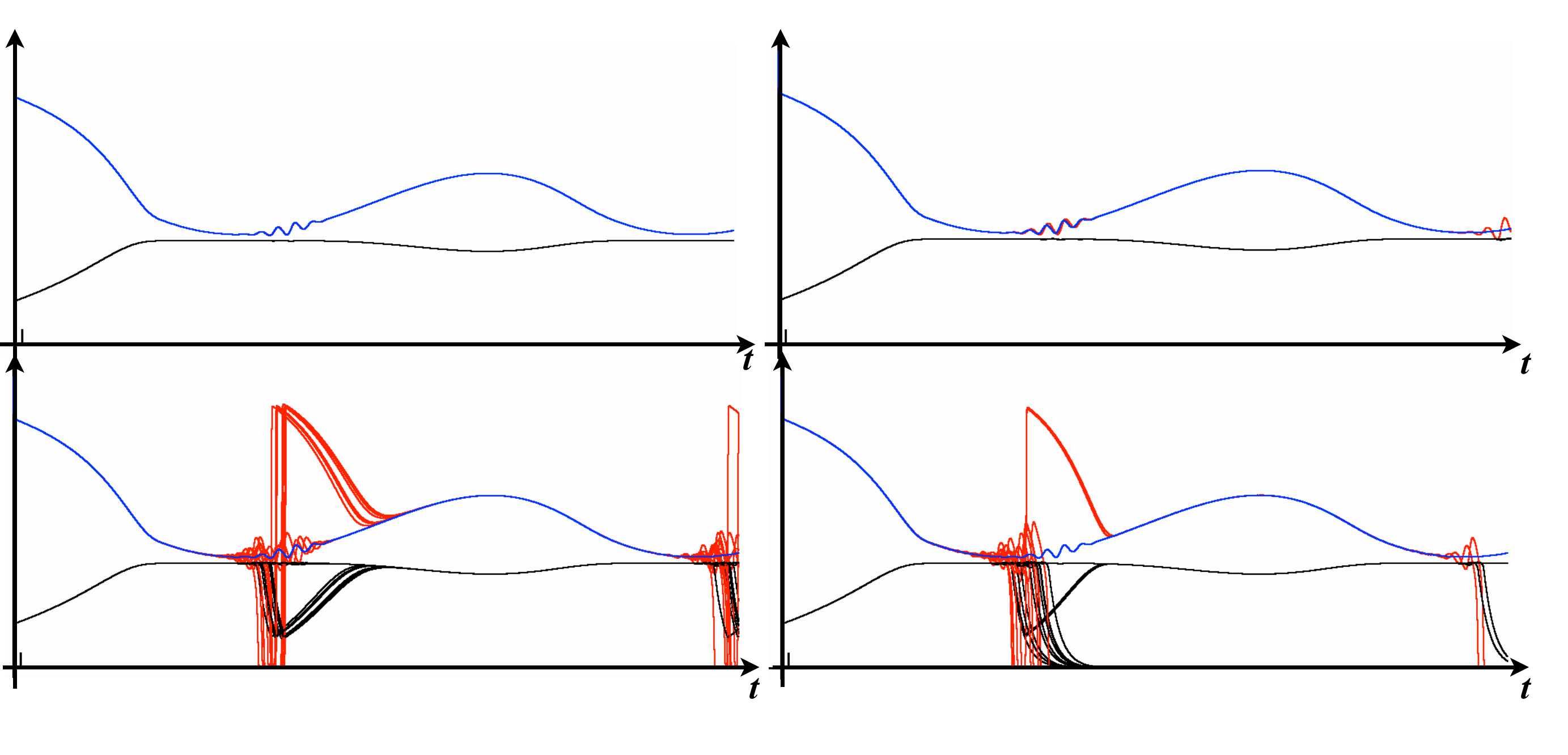} 
   \caption{$\omega = 10^{12}\;10^9$ (above left, right) $\omega = 10^{8}\;10^7$ (below left, right)}
   \label{Fig21}
\end{figure}
An easy way to understand it is to imagine that our parameter $m$ is of the form :
$$ m(t) = a + b\;\mathrm{cos}(r\,t)$$
which mimics, for instance, some seasonal dependence of the mortality rate. This non autonomous system can be considered as a three dimensional system and we see that the ``canard'' value for $m$ is crossed periodically. We have done a simulation in the case :
$$ m(t) = 0.6645 -0.047(1-cos(0.1\,t))$$
and the results are shown on Fig \ref{Fig21}. For $\omega = 10^{12}$ we observe no difference between the deterministic model ($x$ in blue, $y$ in black) and the diffusion approximation ($x$ in red, $y$ in black) ; for $\omega =10^9$ we observe a very slight deviation between red and blue curves ; for $\omega = 10^{8} $ we observe a very big difference with now a mixed mode oscillation in the diffusion process ; for $\omega = 10^7$ the mixed mode oscillation leads to extinction.
\subsection{About the inadequacy of deterministic models with continuous variables.}
\noindent In population dynamics every body agrees that deterministic models are just crude approximations of reality. Only individually based models, stochastic by essence, can represent correctly the evolution of real ecosystems. The example presented here is just one more argument against the danger of using deterministic differential equations without care.\\\\
 {\em But it is by no mean an argument against the study of continuous deterministic differential models of populations dynamics !} \\\\
 Actually there are many good reasons for continuing to explore systems of ordinary differential equations :
 \begin{itemize}
 \item Some models are mathematically appealing. For instance the proof of the exclusion principle for the most general model of competition in the chemostat \cite{SW95}, despite its poor ecological contents, remains an interesting mathematical challenge for mathematicians.
 \item More interesting is the use of easily tractable mathematical models to formalize some ecological issue and clarify the discussion. An interesting example of this use of differential equations is given by the discussion on ``ratio dependent'' models  initiated by the paper of Ardity and Ginsburg \cite{Ard89}. 
 \item In our example, the understanding of the diffusion model, relies on very particular and recently (see bibliographical comments) discovered properties of {\em deterministic} differential systems : the ``canard solutions''.
 \end{itemize}
By the way, far from being an article of propaganda against the use of deterministic differential systems, our paper supports the importance of a thorough understanding of the properties of ordinary differential systems in population dynamics. In particular it shows that the classical deterministic definition of persistence :
$$ \limsup x(t) = \alpha > 0$$
must be enriched  by some consideration about the ``size'' of $\alpha$.
\subsection{About computer simulations in dynamic population modeling.}
\noindent There is no doubt that our mathematical understanding of the phenomena outlined in the present paper will considerably increase in the future. But this will require high mathematical sophistication and time. Unfortunately, in the mean time, biologist will use models and computer simulations which are not completely safe. It urges to provide them with computer routines which are safe of numerical artifacts associated to the true nature of a population : {\em a more or less large number of individuals}. Considering our present mathematical knowledge this certainly can be done in a comparatively short time but it needs quite a lot of people working on the design of safe computer software. This was done in the past for the needs of industry (for instance digital wind tunnels), medicin (medical imaging) this could be the case for microbial ecology but it depends of decisions at the level of scientific policies.
\section{Bibliographical comments.}
\subsection{The atto-fox problem.}
\noindent The question of the inadequacy of deterministic continuous modeling is firmly addressed by D. Mollison \cite{Mol91} in a paper which criticize the biological interpretations of a previous paper by Murray et al. \cite{Mu86}. Let us quote from \cite{Mol91} :
  \begin{quote}
  As to the second wave, close inspection shows that the explanation lies, not much in the determinism of the model, as in its modeling of the population as continuous rather than discrete and its associated inability to let the population variables reach the value zero. Thus the density of infected at the place of origin of the epidemic never becomes zero, it only declines to a minimum of around one atto-fox ($10^{-18}$ of a fox, Hugues 1960) per square kilometer. The model then allows this atto-fox to start the second wave as soon as the susceptible population has regrown sufficiently.
  \end{quote}
  About ten years before Mollison, independently, within the framework of chemical kinetics, D. Gillespie published a famous paper  \cite{GI77} : {\em Exact Stochastic Simulations of Coupled Differential Reactions} from which our model in the present paper is inspired. \\\\
  It is a bit surprising that, at least to our knowledge, not much has been done in this direction. The present paper is a development of a first draft \cite{LoSa09} with T. Sari where we noticed the importance of the presence of ``canards solutions'' regarding the question of persistence in ecological models.  The paper \cite{LORASA09} is also related to this atto-fox question in the case of the chemostat with a slow varying flow rate.  The paper \cite{Cam11} which is much more mathematically oriented, considers the stochastic modeling of the chemostat ; it focusses on the the approximation of jump processes by diffusion processes and was a source of inspiration for the present paper.
  \subsection{Singular perturbations and ``canard solutions''.}
  \noindent As already said, ``canards'' are specific solutions in singular perturbations of differential equations. They where discovered in 1981 by
  a group students of G. Reeb : E. Benoit, J-L. Callot, F. and M. Diener \cite{Ben81}. They studied them within the framework of Non Standard Analysis which is most suitable for modeling since it is a simple formal language where the use of infinitesimals (in the sense where physicists use this term) is mathematically rigorous. But they are now also studied by numerous mathematicians within the framework of {\em mached asymptotic expansion}  or the {\em geometric singular perturbation theory}. The article \cite{Wesh07} by Martin Wechselberger is a short and nice introduction to ``canards'' and the paper \cite{Der10} is a thorough survey about our present understanding of ``canards'' with a focus on numerical questions. The paper \cite{LoSa07} is about Nonstandard Analysis applied to real word questions.\\\\
The question of considering the presence of noise in singularly perturbed systems has been considered for long time. We refer to the recent paper \cite{Ber11}  devoted to the question of the consequence of noisy environment on ``canard solutions'' and its bibliography. In particular the results contained in this paper allow to give asymptotic evaluations of the wide of the ``safety funnel'' and many other quantity of interest but their mathematical sophistication is out of the scope of the present paper.

\section{Conclusion}
\noindent Scientists are now much familiar with the phenomenon of ``sensitivity to initial conditions'' which, in some deterministic dynamical systems, is the cause of an impredictable long range behavior.  The same phenomenon in some deterministic differential equations modeling the dynamic behavior of populations is the cause that a very small difference in an initial condition (or along a trajectory) will make the future value of some variable very small or not. This is the reason why, in the modeling of population dynamic, it is a good thing to add some small noise to the deterministic process because it does not cost too much computer time and may detect this kind of phenomenon. But we have shown that the result may depend strongly of the strength of the noise. By the way, when we do not have an accurate estimation of the strength of the noise, it should be more secure to vary that strength and make sure that the behavior is not strongly dependent on it.

\section{Appendix}
\appendix
\section{Approximation by a diffusion process}

\noindent Consider the process defined by (\ref{lambda}), (\ref{proba}), (\ref{vary}). Since $Z$ follows an exponential law of parameter $\lambda$ its expectation is $\frac{1}{\lambda}$ and the number $N_b$ of events during the duration $dt$ is approximately :
$$ N_b \approx dt /(\frac{1}{\lambda}) = dt \,\lambda = dt \frac{\omega}{\varepsilon}(f(x)+\mu(x)y)$$
We consider $N_b$ as deterministic.
If $dt$ is small the variables $x(t)$ and $y(t)$ are approximatively constant. Denote by $X_i$ the random variable which is equal to one if at the i-th event a predation occurs ; one has :
\begin{equation}
\left \{
 \begin{array}{lcl}

P(X_i = +1) &=&\displaystyle  \frac{\mu(x(t))\,y(t)}{f(x(t))+\mu(x(t))\,y(t)}\\[8pt]
P(X_i = \;\;\;0)& = &\displaystyle \frac{f(x(t))}{f(x(t))+\mu(x(t))\,y(t) }

 \end{array} 
 \right .
 \end{equation}
The number of predations during $[t,t+dt]$ is, approximately, $\sum_1^{N_b}X_i$ and the number of birth is by the way $N_b - \sum_1^{N_b}X_i$ and the increment of the number of individuals is $N_b - 2\sum_1^{N_b}X_i$. \\\\
One has :

\begin{itemize}
\item $\displaystyle  E[\;X_i\;] = \frac{ \mu(x(t))\,y(t)}{f(x(t))+\mu(x(t))\,y(t)}$
\item $\displaystyle E[\;\sum_1^{N_b}X_i\;] = dt \frac{\omega}{\varepsilon}(f(x(t))+\mu(x(t))y(t))\;\;\frac{ \mu(x(t))\,y(t)}{f(x(t))+\mu(x(t))\,y(t)} = dt \frac{\omega}{\varepsilon} \;\mu(x(t))\,y(t)$
\item $\displaystyle \sigma^2(X_i) = \frac{f(x(t))\mu(x(t))\,y(t)}{(f(x(t))+\mu(x(t))\,y(t))^2}          $
\item $\displaystyle \sigma^2(\;\sum_1^{N_b}X_i\;) = dt \frac{\omega}{\varepsilon}(f(x(t))+\mu(x(t))y(t)) \;\;\;\frac{f(x(t))\mu(x(t))\,y(t)}{(f(x(t))+\mu(x(t))\,y(t))^2} $
\item$\displaystyle \sigma^2(\;\sum_1^{N_b}X_i\;) = dt \frac{\omega}{\varepsilon}\;\;\;\frac{f(x(t))\mu(x(t))\,y(t)}{(f(x(t))+\mu(x(t))\,y(t))}$
\end{itemize}
From the central limit theorem we can approximate the sum by a Gaussian and we write :
$$ \sum_1^{N_b}X_i\ \approx dt \frac{\omega}{\varepsilon} \;\mu(x(t))\,y(t) + \sqrt{dt \frac{\omega}{\varepsilon}\;\;\;\frac{f(x(t))\mu(x(t))\,y(t)}{(f(x(t))+\mu(x(t))\,y(t))}}W_t$$
where $W_t$ is a Gaussian of $0$ mean and  $1$ as standard deviation . \\\\
Since the variable $x$ is the number of individuals divided by $\omega$ the increment of $x$ is given by :
$$x(t+dt) - x(t) = \frac{1}{\omega} (N_b - 2\sum_1^{N_b}X_i)$$
$$\displaystyle x(t+dt) - x(t) \approx \frac{1}{\omega} \left \{ N_b - 2\left \{dt \frac{\omega}{\varepsilon} \;\mu(x(t))\,y(t) + \sqrt{dt \frac{\omega}{\varepsilon}\;\;\;\frac{f(x(t))\mu(x(t))\,y(t)}{(f(x(t))+\mu(x(t))\,y(t))}}W_t \right \}\right \}$$
and replacing by the value of $N_b$ one get :
\begin{figure}[t,b]   \centering
   \includegraphics[width=12cm]{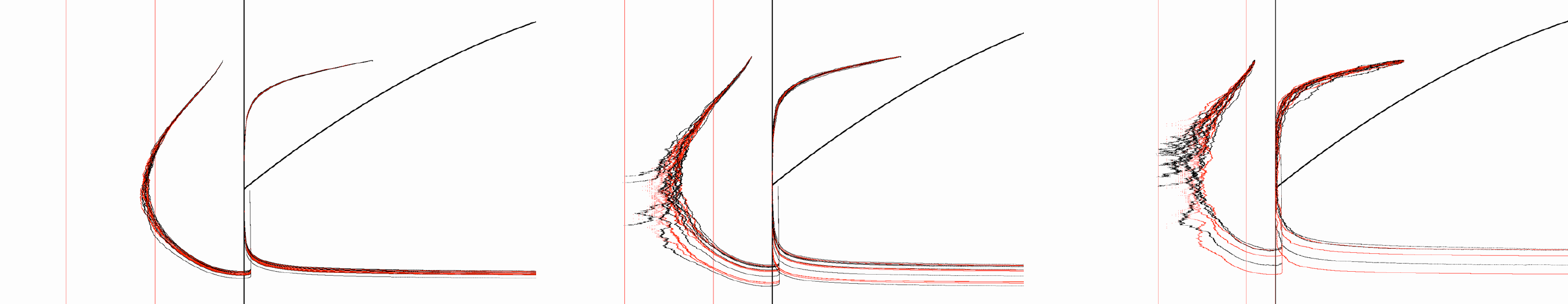} 
   \caption{Comparison between birth and death process and diffusion process }
   \label{D-G}
\end{figure}

$$\displaystyle x(t+dt) - x(t) \approx dt \frac{1}{\varepsilon} [f(x(t))-\mu(x(t))\,y(t)] - \sqrt{dt \frac{4}{\omega \varepsilon}\;\;\;\frac{f(x(t))\mu(x(t))\,y(t)}{(f(x(t))+\mu(x(t))\,y(t))}}W_t $$
Let us compute now the increment of  $y$.
According to (\ref{vary}) we have :
$$y(t+dt) = y(t) - dt\,m\, y(t) + \varepsilon \{ \mathrm{number \;of\;prey\;death\; during} \; [t, t+dt] \}$$
which, according to the previous notations is :
$$y(t+dt) - y(t)  = -  dt\,m\, y(t) + \frac{\varepsilon}{\omega} \sum_1^{N_b}X_i$$
and introducing $W_t$ one gets :
$$y(t+dt) - y(t)  \approx    dt[ \mu(x(t))-m] y(t) + \sqrt{dt \frac{\varepsilon}{\omega}\;\;\;\frac{f(x(t))\mu(x(t))\,y(t)}{(f(x(t))+\mu(x(t))\,y(t))}}W_t $$
On Fig.\ref{D-G} one sees a comparison between the birth and death process (red trajectories) and its approximation by a diffusion. From the left to the right we have $\omega = 10^6$, $\omega = 10^5$, $\omega = 10^4$. The representation is both in $(x,y)$ and $(\xi,y)$ variables. We have $10$ runs from the initial condition $(0.2\,,\, 0.6)$. The two red vertical lines correspond to a population between $1$ and $1000$.

\section{Exponentially small values}
\noindent Let us write explicitly system (\ref{syst4}) as :

\begin{equation} \label{explicit}
\left \{
 \begin{array}{lcl}
 \displaystyle \frac{dx}{dt}& = &\displaystyle \frac{1}{\varepsilon}[0.5 x(2-x) -\frac{x}{0.4+x}y ] \\[8pt]
 \displaystyle \frac{dy}{dt}&=&\displaystyle  (\frac{x}{0.4+x} - m)y
 \end{array} 
\right .
\end{equation}
In the variables $(\xi,y)$ the system writes :
$$
\left \{
 \begin{array}{lcl}
 \displaystyle \frac{d\xi}{dt}& = &\displaystyle [0.5(2-(\xi/\varepsilon)) -\frac{1}{0.4+(\xi/\varepsilon)}y ] \\[8pt]
 \displaystyle \frac{dy}{dt}&=&\displaystyle  (\frac{\xi/\varepsilon)}{0.4+\xi/\varepsilon)} - my
 \end{array} 
\right .
$$
which is approximated, when $\xi << \varepsilon$, by : 
\begin{equation} \label{logvar}
\left \{
 \begin{array}{lcl}
 \displaystyle \frac{d\xi}{dt}& = &[1-2.5 y ] \\[8pt]
 \displaystyle \frac{dy}{dt}&=& - my
 \end{array} 
\right .
\end{equation}
Take as initial condition $( x_0, y_0) = (-0.1,0.9)$ (which corresponds to trajectory $n^o \;1$ in Fig. \ref{log1}) and integrate. It comes that the minimum $\xi^*$ for $\xi(t)$ is attained for the value $t^*$ of $t$ for which $y(t^*) = 0.4$ and this value turns out to be approximately $-1$. But :
$$ x^*  =e^{ \frac{\xi^*}{0.02} } \approx e^{-40}\approx 10^{-17}$$
The minimum depends much of the value of  $y_0$ : The largest is $y_0$ the smallest is the minimum. This explain why in Fig.\ref{log2} the minimum corresponding to trajectory $6$ is much bigger.

\section{Numerical simulations}
\noindent We did not use any solver. A specific software was written in order to be sure that there were not artifacts caused by erroneous uses of some sophisticated numerical scheme. Trajectories of the differential equations (\ref{syst4}) are obtained using the Euler scheme defined by (\ref{syst3}). We prefer this scheme to any more sophisticated scheme used to simulate differential systems since it is the exact recurrence scheme which approximate for  $E[x(t)]$ and $E[y(t)]$ of the diffusion process (\ref{syst2}).\\\\
We fixed $dt = 10^{-4}$ since we observed that for this value solutions of  (\ref{syst3}) are indistinguishable from those with $dt = 10^{-5}$.\\\\
The birth and death process defined by (\ref{lambda}), (\ref{proba}), (\ref{vary}) takes too long time  to be simulated when $\lambda$ is very large
(in the case of our computer $\omega > 10^{6}$) and this is the reason why we used a diffusion approximation which is a perfect approximation for large values. Since we where mainly interested in the funnel phenomenon associated to ``canard '' it was not necessary to switch to the true birth and death process for small values of $\lambda$. But if one is interested by figures like the mean of the extinction time it should be better to switch to some suitable jump process.

\subsection*{Acknowledgements.}
\noindent The ``deterministic part'' of the paper is inspired by the paper \cite{LoSa09} of the second author with T. Sari. Since the publication of this paper we had many fruitful discussions with him and also with J. Harmand and A. Rapaport from the Modemic Team (http://www-sop.inria.fr/modemic/). We thanks them warmly.\\\\
The financial support of the French National Research Agency (ANR) within the SYSCOMM project DISCO ANR-09-SYSC-003.
 is appreciated.

\end{document}